\documentclass{elsart3-1}
\usepackage{amsfonts}
\usepackage{amsmath}
\usepackage{amssymb}
\usepackage[utf8]{inputenc}
\usepackage[english,french]{babel}  
\usepackage[T1]{fontenc}      

\usepackage{color}

\def\red#1{\textcolor{red}{#1}} 
\def\red#1{#1}

\newtheorem{e-proposition}[theorem]{Proposition}

\newtheorem{e-definition}[theorem]{Definition\rm}

\newtheorem{theoreme}{Th\'eor\`eme}[section]

\newtheorem{remarque}{\it Remarque}

\def\og{\leavevmode\raise.3ex\hbox{$\scriptscriptstyle\langle\!\langle$~}}
\def\fg{\leavevmode\raise.3ex\hbox{~$\!\scriptscriptstyle\,\rangle\!\rangle$}}

\def\Xint#1{\mathchoice
  {\XXint\displaystyle\textstyle{#1}}%
  {\XXint\textstyle\scriptstyle{#1}}%
  {\XXint\scriptstyle\scriptscriptstyle{#1}}%
  {\XXint\scriptscriptstyle\scriptscriptstyle{#1}}%
  \!\int}
\def\XXint#1#2#3{{\setbox0=\hbox{$#1{#2#3}{\int}$}
    \vcenter{\hbox{$#2#3$}}\kern-.5\wd0}}

\def\fint{\Xint-}

\journal{the Acad\'emie des sciences}
\begin{document}
\centerline{}
\begin{frontmatter}




%
\selectlanguage{french}
\title{Approximation locale précisée dans des problèmes multi-échelles avec défauts localisés}



\author[authorlabel1]{Xavier Blanc},
\ead{blanc@ann.jussieu.fr}
\author[authorlabel2]{Marc Josien},
\ead{marc.josien@enpc.fr}
\author[authorlabel2]{Claude Le Bris},
\ead{claude.le-bris@enpc.fr}

\address[authorlabel1]{Univ. Paris Diderot, Sorbonne Paris Cité, Laboratoire Jacques-Louis Lions, UMR 7598, UPMC, CNRS, F-75205 Paris, FRANCE}
\address[authorlabel2]{Ecole des Ponts and INRIA, 6 \& 8, avenue Blaise Pascal, 77455 Marne-La-Vall\'ee
  Cedex 2, FRANCE}


\medskip
\selectlanguage{french}
\begin{center}
{\small Re\c{c}u le *****~; accept\'e apr\`es r\'evision le +++++\\
Pr\'esent\'e par £££££}
\end{center}

\begin{abstract}
\selectlanguage{french}
Nous poursuivons l'étude initiée dans~\cite{BLLMilan} de problèmes multi-échelles avec défauts, dans
le cadre de la théorie de l'homogénéisation, spécifiquement ici pour une équation de diffusion  avec un coefficient de la forme fonction périodique perturbée par une fonction $L^r (\mathbb{R}^d)$, $1<r< +\infty$, modélisant un défaut local. Nous esquissons la démonstration du fait  que le correcteur, dont l'existence a été prouvée dans~\cite{BLLMilan,BLLcpde}, permet d'approcher la fonction solution de l'équation originale avec la même précision, essentiellement, que dans le cas purement périodique. Les taux de convergence varient, et sont précisés, en fonction de l'intégrabilité $L^r$ du défaut. Une extension à un cas abstrait "général" est mentionnée. Les résultats annoncés dans cette Note seront précisés dans les documents~\cite{Article_B,Article_A}.
\vskip 0.5\baselineskip

\selectlanguage{english}
\noindent{\bf Abstract}
\vskip 0.5\baselineskip
\noindent
{\bf Local precised approximation in multiscale problems with local defects.}

We proceed here with our systematic study, initiated in~\cite{BLLMilan}, of multiscale problems with defects,
within the context of homogenization theory. The case under consideration here is that of a diffusion equation with
a diffusion coefficient of the form of a periodic function perturbed by an $L^r (\mathbb{R}^d)$ , $1<r< +\infty$,
function modeling a localized defect. We outline the proof of the following approximation result: the corrector
function, the existence of which has been established in~\cite{BLLMilan,BLLcpde}, allows to approximate the
solution of the original multiscale equation with essentially the same accuracy as in the purely periodic case. The
rates of convergence may however vary, and are made precise, depending upon the  $L^r$ integrability of the
defect. The generalization to an abstract setting is mentioned. Our proof exactly follows, step by step, the pattern of the original proof of
Avellaneda and Lin in~\cite{AvellanedaLin} in the periodic case, extended in the works of Kenig and
collaborators~\cite{KLSGreenNeumann}, and borrows a lot from it. The details of the results announced in this Note are given in our publications~\cite{Article_B,Article_A}.
\end{abstract}
\end{frontmatter}

\selectlanguage{english}
\section*{Abridged English version}
We continue our study \cite{BLLMilan,BLLcpde} of homogenization problems in nonperiodic media. The particular setting considered here is
that of a diffusion equation~\eqref{DefUepsilon} with a  diffusion coefficient $a$ of the
form~\eqref{eq:forme_coefficient}, where $\widetilde a$ is an $L^r (\mathbb R^d)$ , $1<r< +\infty$, function modeling a
localized defect, decaying at infinity in a loose sense, and perturbing the background periodic medium $a_{\rm
  per}$. We aim at quantitatively estimating  at which rate the two-scale expansion (truncated at the first order)
provided by homogenization theory approaches the exact solution $u^\epsilon$ as $\epsilon$ vanishes. Put
differently, we seek the rate at which the remainder term~\eqref{DefResteEps} goes to zero. In~\eqref{DefResteEps},
the corrector employed is the function $w_p$ solution to the corrector equation~\eqref{ProbCorrPer}, the existence
and uniqueness (up to additive constants) of which has been established in our previous
works~\cite{BLLMilan,BLLcpde}. Such a corrector is different from the periodic corrector, and although intuitively
one could have thought, based on the observation that the presence of $\widetilde a$ does not modify the homogenized
equation~\eqref{DefUetoile}, that the periodic corrector $w^{per}_p$ would give an equally accurate approximation, it does
not. The rates of convergence obtained are made precise in our main result, namely Théorème 2.1 of the French
version, and estimates~\eqref{PremRes} through~\eqref{TroisRes} in various norms. Interestingly, the rates of
convergence may however vary from one case to another, and also in comparison with the periodic case. They depend
upon the  $L^r$ integrability of the defect. Our proof exactly follows, step by step, the pattern of the original proof of
Avellaneda and Lin in~\cite{AvellanedaLin} in the periodic case, extended in the works of Kenig and
collaborators~\cite{KLSGreenNeumann}, and borrows a lot from it. Instead of having a bounded (periodic) corrector,
as in those proofs, we have here a corrector function that is not necessarily bounded. The crucial ingredient of
the proof is then, in fact,  the strict sublinearity of the corrector function, which, given our assumptions, can
be made precise, see~\eqref{CorrSlin}. This suggests a generalization of our setting, beyond the "periodic +local
perturbation" case,  which we make precise in Section 3.1 of the French version: essentially, besides usual
assumptions, the key point is that the corrector is stricly sublinear at infinity with a prescribed rate of
sublinearity and (a property that is actually very much linked to the former) that the potential function
associated to this corrector (defined in~\eqref{AssocAPot}-\eqref{eq:potentiel}) is also strictly sublinear at
infinity in a similar quantitative manner (see~\eqref{PotSlin}).  In passing, and as is also the case in the proof
of the periodic case, we establish estimates for the Green function $G^\epsilon$ of the original problem
(see~\eqref{Green1Dir}, \eqref{Green2Dir}, \eqref{Green3Dir} of the French version) and its convergence to the (possibily corrected) Green function of the homogenized problem~\eqref{DefUetoile} (see~\eqref{DefGreen1}). The details of the results announced in this Note are given in our publications~\cite{Article_B,Article_A}.

\selectlanguage{french}

\newpage
\section{Introduction}

  \subsection{Motivation}

    Dans cette Note, nous poursuivons l'étude initiée dans \cite{BLLMilan} de problèmes elliptiques multi-échelles,
    dans le cadre de la théorie de l'homogénéisation. L'équation que nous considérons possède un coefficient qui
    présente, à l'échelle microscopique, une structure périodique perturbée localement par un \og défaut \fg. On se
    convainc aisément que le comportement macroscopique d'un tel matériau est dicté par sa seule structure périodique
    sous-jacente. Si, en revanche, on cherche à obtenir une information plus fine, en terme de taux de convergence, à l'échelle
    microscopique, pour une norme plus forte ou encore au voisinage du défaut, alors ce défaut ne peut plus être
    négligé. 
    
Dans \cite{BLLMilan}, il a été montré (en 1D au moins, et formellement en dimension supérieure) dans un cadre hilbertien, i.e pour un défaut dans $L^2(\mathbb{R}^d)$, qu'il
est en effet nécessaire de construire un correcteur prenant en compte le défaut pour obtenir une approximation
précisée de la solution. L'existence d'un tel correcteur est démontrée, et ce résultat est généralisé dans
\cite{BLLcpde} au cas d'un défaut d'intégrabilité $L^r(\mathbb{R}^d)$, $1<r<+\infty$, ainsi qu'à d'autres
situations "perturbatives". Il est affirmé dans
\cite{BLLMilan,BLLcpde} que, formellement, un tel correcteur permet d'assurer l'approximation voulue. L'objet de
cette Note est d'énoncer précisément ce résultat, et de donner les grandes lignes de sa preuve : \og le correcteur
corrige\fg{} dans la norme considérée, à un ordre précisé. Les résultats annoncés dans cette Note, et leurs preuves, seront détaillés dans les publications \cite{Article_B,Article_A} en préparation.

  \subsection{Le cas périodique}
    Nous nous donnons un champ $a\in {L}^\infty\left(\mathbb{R}^d\right)$, pris à valeurs scalaires pour simplifier l'exposé,  et nous considérons le problème suivant:
    \begin{align}
      \label{DefUepsilon}
      -\operatorname{div} (a(x/\epsilon)\nabla u^\epsilon(x)) =f(x) \quad\text{dans}\quad \Omega, \quad\text{et}\quad u^\epsilon=0 \quad\text{sur}\quad \partial \Omega, 
    \end{align}
    posé sur un domaine borné régulier $\Omega \subset \mathbb{R}^d$. Le champ $a$ est elliptique. Il est bien connu (voir par exemple \cite{JKO}) que, dans le cas où $a$ est périodique, alors le problème \eqref{DefUepsilon} s'homogénéise en le problème suivant:
    \begin{align}
      \label{DefUetoile}
      -\operatorname{div}\left(a^* . \nabla u^*(x)\right)  =f(x)  \quad\text{dans}\quad \Omega, \quad\text{et}\quad u^*=0 \quad\text{sur}\quad \partial \Omega,
    \end{align}
    où $a^*$ est une matrice constante. En particulier, on observe la convergence faible $\nabla u^\epsilon\rightharpoonup \nabla u^*$ dans ${L}^2(\Omega)$. Pour obtenir de la convergence forte, il faut corriger $u^*$.
    Pour ce faire, on définit les \textit{correcteurs} $w_j$, $j \in [\![1,d]\!]$, associés à $a$ comme étant les solutions de l'équation suivante:
     \begin{align}
      \label{ProbCorr}
	  -\operatorname{div}\left(a\left( e_j+ \nabla w_j \right) \right)=0 & \quad\text{dans}\quad \mathbb{R}^d,  \quad\text{et}\quad |w_j(x)|/(1+|x|) \underset{|x| \rightarrow + \infty}{\rightarrow} 0,
    \end{align}
    où les $e_j$ sont les vecteurs de la base canonique de $\mathbb{R}^d$, et on introduit le reste défini par
    \begin{align}\label{DefResteEps}
	R^\epsilon(x):=u^\epsilon(x)-u^*(x) - \epsilon \sum_{j=1}^d w_j\left( \frac{x}{\epsilon} \right) \partial_j u^*(x).
    \end{align}
Les correcteurs $w_j$ sont périodiques, et on démontre classiquement que $\left\|\nabla R^\epsilon \right\|_{{L}^2(\Omega)} \rightarrow 0$. 
    Avellaneda et Lin ont aussi démontré dans \cite{AvellanedaLin} que, si de plus $a$ est Hölderienne, on peut
    obtenir des estimations lipschitziennes sur $u^\epsilon$, et des estimations sur $R^\epsilon$ quantifiées en
    $\epsilon$. Des estimations similaires ont été obtenues dans le cas de coefficients non réguliers stationnaires dans \cite[Corollary 3]{gloria}.

Ces travaux ont notamment été approfondis dans \cite{KLSGreenNeumann}, où est démontré que, modulo le fait de prendre des correcteurs $w_j$ adaptés au domaine (c'est à dire avec une définition légèrement différente de \eqref{ProbCorr}), on obtient l'estimation suivante: $\left\|\nabla R^\epsilon\right\|_{{L}^\infty(\Omega)}\leq C\epsilon \ln \epsilon$ (voir \cite[Lem.\ 3.5]{KLSGreenNeumann}). Au cours de la preuve, les auteurs approximent la fonction de Green $G^\epsilon$ associée au problème \eqref{DefUepsilon}, ses gradients $\nabla_xG^\epsilon$ et $\nabla_y G^\epsilon$, et son gradient croisé $\nabla_x\nabla_y G^\epsilon$. Pour ce faire, ils emploient la fonction de Green $G^*$ du problème \eqref{DefUetoile}, convenablement modifiée par les correcteurs $w_j$ (voir \cite[Th.\ 3.6 et Th. 3.11]{KLSGreenNeumann}).

  \subsection{Le cas périodique perturbé par un défaut local}
  
    Dans \cite{BLLMilan,BLLcpde} est considéré le cas d'un champ scalaire 
    \begin{align}\label{eq:forme_coefficient}
      a=a_{\rm per}+\widetilde{a},
    \end{align}
    où $a_{\rm per}$ est périodique et $\widetilde{a}$ est une perturbation locale. 
 Plus précisément, supposons que $d \geq 3$ (en dimension $2$, des détails techniques supplémentaires sont
 nécessaires du fait que les fonctions de Green d'opérateurs elliptiques ne sont pas bornées à l'infini) et qu'il existe $\alpha>0$, $\mu>0$ et $r\in ]1,+\infty[$ tels que
    \begin{align}
      \label{HypoBLL}
      a_{\rm per}\ , \ \widetilde{a} \in C_{\rm unif}^{0,\alpha}\left(\mathbb{R}^d\right),\quad \mu^{-1}
      \leq a_{\rm per} \leq \mu, \quad \mu^{-1} \leq \widetilde{a} + a_{\rm per}\leq \mu, \quad \text{ et }\quad \widetilde{a} \in {L}^r\left(\mathbb{R}^d\right).
    \end{align}
Ici, $C_{\rm unif}^{0,\alpha}\left(\mathbb{R}^d\right)$ désigne l'espace des fonctions uniformément
Hölderiennes sur ${\mathbb R}^d$, de coefficient $\alpha\in ]0,1[$. 
    Alors par \cite[Th.\ 4.1]{BLLcpde}, il existe une solution $w_j$ au probl\`eme \eqref{ProbCorr} qui s'\'ecrit $w_j=w^{\rm per}_j + \widetilde{w}_j$, où $\nabla \widetilde{w}_j \in {L}^r\left(\mathbb{R}^d\right) \cap {L}^\infty\left(\mathbb{R}^d\right)$ et où $w^{\rm per}_j$ est une solution p\'eriodique de
    \begin{equation}
      \label{ProbCorrPer}
	  -\operatorname{div}\left(a_{\rm per}\left( e_j+\nabla w^{\rm per}_j \right) \right)=0 \quad\text{dans}\quad \mathbb{R}^d.
    \end{equation}
    Par le Théorème de Morrey \red{\cite[Th.~9.12]{brezis},} cela implique en particulier que, si $r\neq d$, les correcteurs $w_j$ satisfont
    \begin{equation}\label{CorrSlin}
      \left|w_j(x)-w_j(y) \right| \leq C \left|x-y\right|^{1-{\nu}}, 
    \end{equation}
    pour tous $x, y \in \mathbb{R}^d$, où
    \begin{align}
      \label{DefNur}
      \nu = {\nu_r}:=\min\left(1 ,d/r \right) \in ]0,1].
    \end{align}
Dans le cas $r=d$, qui est critique, on obtient \eqref{CorrSlin} seulement pour tout $\nu<1$. Dans tous les cas,
cette estimation est \og seulement contraignante\fg{} pour $|x-y|$ grand. En effet, pour $|x-y|$ petit,  $w_j\in
C^{1,\alpha}_{\rm unif}$ par régularité elliptique (pour la valeur $\alpha$ de \eqref{HypoBLL}), et \eqref{CorrSlin} est
triviale (avec $\nu=0$). L'inégalité \eqref{CorrSlin} traduit une \og sous-linéarité renforcée\fg~ des correcteurs (qui devient, pour $r<d$, une borne ${L}^\infty$ sur lesdits correcteurs).
    Comme nous allons le voir, l'existence d'un correcteur ainsi que l'estimation \eqref{CorrSlin} sont des ingrédients essentiels pour contrôler le gradient du reste $\nabla R^\epsilon$.

\section{Résultats}

    Nous démontrons un résultat qui étend au cas ci-dessus une partie de l'analyse faite dans \cite{AvellanedaLin}
    et \cite{KLSGreenNeumann} (en particulier le Théorème~5 de \cite{AvellanedaLin} et les Théorèmes~3.4 et 3.7 de
    \cite{KLSGreenNeumann}) :
    \begin{theoreme}\label{ThDefautLp}
      Soit $d\geq 3$.
      Supposons qu'il existe $\alpha>0$, $\mu>0$ et $r\in ]1,+\infty[$, $r\neq d$ tels que $a_{\rm per}$  et  $\widetilde{a}$ satisfont \eqref{HypoBLL}, et soit ${\nu_r}$ défini par \eqref{DefNur}.
      Soient $\Omega$ un domaine régulier et $\Omega_1 \subset \subset \Omega$.
      Considérons $f \in {L}^2(\Omega)$ et $u^\epsilon, u^*, R^\epsilon$ respectivement définies par \eqref{DefUepsilon}, \eqref{DefUetoile}, et \eqref{DefResteEps}.
      Alors $R^\epsilon \in {H}^1(\Omega)$ et
      \begin{align}
	\label{PremRes}
	&\left\|R^\epsilon\right\|_{{L}^2(\Omega)} \leq C_1 \epsilon^{{\nu_r}} \left\|f\right\|_{{L}^2(\Omega)},
	\\
	\label{DeuxRes}
	&\left\|\nabla R^\epsilon\right\|_{{L}^2(\Omega_1)^d} \leq C_2 \epsilon^{{\nu_r}} \left\|f\right\|_{{L}^2(\Omega)}.
      \end{align}
      En outre, pour tout $\beta \in ]0,\alpha]$, si $f \in {C}^{0,\beta}\left(\bar{\Omega}\right)$, on a $R^\epsilon \in {W}^{1,\infty}(\Omega)$ et
      \begin{equation}
	\label{TroisRes}
	\left\| \nabla R^\epsilon \right\|_{{L}^\infty(\Omega_1)^d} \leq C_3 \epsilon^{{\nu_r}}\ln\left(2+\epsilon^{-1}\right)\left\|f\right\|_{{C}^{0,\beta}(\Omega)^d},
      \end{equation}
où les constantes $C_1,C_2,C_3$ sont indépendantes de $\epsilon$ et $f$. 
    \end{theoreme}
L'intégrabilité $L^r$ du défaut détermine la qualité de l'approximation. Le cas $r=1$ n'est pas inclus dans
l'énoncé ci-dessus car on ne sait pas alors démontrer l'existence d'un correcteur dans l'espace fonctionnel
associé. Toutefois, on peut obtenir le résultat en notant que $\tilde a\in L^r$ pour tout $r\geq 1$, et appliquer
le théorème~\ref{ThDefautLp} avec $r<d$, qui donne $\nu_r=1$. D'autre part, $r=d$ est critique, ce qui est naturel (voir \cite[Section 3]{BLLcpde}). Il constitue la charnière entre
deux régimes. En effet, si $r<d$, les résultats sont les mêmes que dans le cadre de l'homogénéisation périodique,
et ce parce que $w_j\in {L}^\infty\left(\mathbb{R}^d\right)$. En revanche, si $r>d$, alors le correcteur n'est
pas borné a priori, d'où un moindre contrôle sur la quantité $R^\epsilon$: le défaut devient suffisamment "gros" pour
avoir un impact macroscopique (pour l'approximation de $u^\epsilon$ dans des normes assez fines et/ou à l'ordre $\epsilon$). 
    Dans le cas $r=d$, il n'est pas prouvé, ni même certain, que les correcteurs $w_j$ sont bornés (voir \cite{BLLcpde}). Mais on peut se ramener (de façon sous-optimale) aux résultats prouvés pour $r>d$, puisque $L^d\cap
    L^\infty\subset L^r\cap L^\infty$ pour $d<r<+\infty$.
    Notons enfin que la présence du correcteur dans \eqref{PremRes} est superflue, et qu'on peut avoir la même estimation sur $u^\epsilon-u^*$ (au lieu de $R^\epsilon$).

\section{Remarques et extensions possibles}

  Nous faisons ici quelques remarques et renvoyons aux publications en préparation \cite{Article_B,Article_A} pour plus de précisions.

\subsection{Cadre abstrait général}
    La démonstration du Théorème~\ref{ThDefautLp} ne fait usage de l'hypothèse de la 
structure "périodique + défaut" que pour démontrer
    l'existence de correcteurs et d'un potentiel (à savoir la fonction $B$ définie en
    \eqref{AssocAPot}-\eqref{eq:potentiel} ci-dessous) fortement
    sous-linéaires. Ainsi, les conclusions du Théorème~\ref{ThDefautLp} sont en fait valides sous les hypothèses
    suivantes, plus générales que celles utilisées ici :
    \begin{enumerate}
    \item la matrice $a$ est elliptique, bornée, uniformément Hölderienne ;
    \item la matrice homogénéisée $a^*$ est constante ; 
    \item elle admet un correcteur $w_j$, c'est-à-dire une solution de \eqref{ProbCorr} ;
    \item ce correcteur est fortement sous-linéaire à l'infini, c'est-à-dire qu'il vérifie \eqref{CorrSlin}, pour
      un certain $\nu\in ]0,1]$ ;
    \item il existe un potentiel $B$ associé (i.e une solution antisymétrique de~\eqref{AssocAPot}-\eqref{eq:potentiel} ci-dessous), qui est lui
      aussi fortement sous-linéaire, c'est-à-dire qu'il vérifie \eqref{PotSlin}, pour $\nu\in]0,1]$.
    \end{enumerate}
Ces hypothèses impliquent en particulier que $a(x/\epsilon)$ H-converge uniformément vers $a^*$, propriété que l'on
définit comme suit : pour tout suite $\epsilon_n\to 0$ et toute suite $\left(y_n\right)_{n\in\mathbb N}$ de
$\mathbb R^d$, 
\begin{equation}
  \label{eq:unifHcv}
 a\left(\frac x {\epsilon_n} + y_n\right)\quad \text{H-converge vers}\quad a^*.  
\end{equation}
 Pour la
définition de la H-convergence, nous renvoyons à \cite[Definition 6.4]{Tartar}. Cette propriété est fondamentale dans la preuve esquissée ci-dessous. Pour les détails
de cette généralisation, nous renvoyons encore à \cite{Article_B,Article_A}.

\subsection{Autres remarques}

\label{sec:remarq-et-extens}

  \begin{enumerate}
    \item{La preuve esquissée ici est faite dans le cas où $a$ est scalaire. Toutefois, il est possible de
        travailler avec un coefficient matriciel, \red{toujours pour une équation scalaire.} On obtient alors des résultats analogues.}
    \item{La preuve originale de \cite{AvellanedaLin} fonctionne pour des systèmes. Par conséquent, dans la mesure
        où l'existence des correcteurs $w_j$ est aussi prouvée dans \cite{BLLfutur1} pour le cas des systèmes, il semble a priori possible de
        démontrer un résultat analogue au Théorème~\ref{ThDefautLp} dans le cadre d'un système d'équations.
        Une telle adaptation n'a cependant pas été entreprise. Voir à ce sujet la Remarque~\ref{rk1} ci-dessous. Notons que, dans le cas d'une équation, le principe du maximum et le Théorème de De Giorgi-Nash-Moser \cite[Th.\ 8.24 p.\ 202]{GT} permettent de simplifier certains aspects de la démonstration.}
    \item{De la même manière que dans \cite{KLSGreenNeumann}, il est possible d'approximer la fonction de Green $G^\epsilon$ relative à l'Equation \eqref{DefUepsilon}, ainsi que ses gradients $\nabla_x G^\epsilon$ et $\nabla_y G^\epsilon$, et son gradient croisé $\nabla_x \nabla_y G^\epsilon$.
    }
    \item{
    Les estimations \eqref{DeuxRes} et \eqref{TroisRes} sont des estimations à l'intérieur du domaine. Toutefois, dans le cas d'une matrice périodique, on obtient des estimations jusqu'au bord (voir \cite{KLSGreenNeumann}). Cela requiert d'introduire des correcteurs adaptés au domaine, lesquels peuvent être construits à partir des correcteurs définis sur tout $\mathbb{R}^d$ (voir \cite[Prop.\ 2.4]{KLSGreenNeumann}). 
    Ces correcteurs sont également bien définis dans le cas présent, ce qui fournit un point de départ pour adapter la preuve de \cite{KLSGreenNeumann}.
    }
    \item{\label{item:1}
    On peut aussi montrer dans le cadre du Théorème~\ref{ThDefautLp} que, si $f\in {L}^p(\Omega)$, pour tout $p\in [2,+\infty[$, alors
    \begin{equation*}
       \left\| R^\epsilon\right\|_{{L}^p(\Omega)} \leq C \epsilon^{\nu_r} \left\| f\right\|_{{L}^p(\Omega)} \quad\text{et}\quad \left\| \nabla R^\epsilon\right\|_{{L}^p(\Omega_1)} \leq C \epsilon^{{\nu_r}} \left\|f\right\|_{{L}^p(\Omega)}.
    \end{equation*}
    L'estimation sur $R^\epsilon$ est immédiate vu le schéma de preuve ci-dessous. 
    L'estimation sur $\nabla R^\epsilon$ découle d'un Lemme de mesure à la Calder\'on-Zygmund (voir \cite[Th.\ 2.4]{ShenCalderon}).}
  \end{enumerate}

\section{Schéma de preuve}

  Notre schéma de preuve suit celui des articles \cite{AvellanedaLin,KLSGreenNeumann}. L'idée repose sur le fait
  que, pour $\epsilon=0$, l'équation est à coefficients constants, donc vérifie des estimations de régularité
  elliptique, à la fois de type Schauder (en normes $C^{k,\alpha}$), et de type Sobolev (en normes
  $W^{k,p}$). Pour $\epsilon$ petit, par compacité, on arrive à obtenir des propriétés similaires. Ceci est l'idée
  fondamentale des preuves de \cite{AvellanedaLin}, laquelle repose sur le fait que le correcteur $w_j$ est borné,
  et donne, via l'expression \eqref{DefResteEps}, une bonne approximation de $u_\epsilon$. Ici, le correcteur n'est
  plus borné a priori, mais le fait que $\nabla \widetilde w_j\in L^r(\mathbb{R}^d)$ nous permet d'adapter les preuves de
  \cite{AvellanedaLin,KLSGreenNeumann}. Notons que le cas où le défaut est d'intégrabilité $r<d$ ne nécessite que des modifications très minimes, car dans ce cas on sait que les $w_j$ sont bornés, ce qui est un ingrédient fondamental des démonstrations de \cite{AvellanedaLin}.
  En revanche, si $r>d$, il faut faire quelques modifications ponctuelles et techniques, qui changent notamment le
  taux d'approximation (d'où la présence de l'exposant ${\nu_r}$ dans le Théorème~\ref{ThDefautLp}). Nous donnons
  ci-dessous les grandes lignes de la démonstration, en indiquant en particulier les points où le fait que $w_j$
  n'est pas borné nécessite des adaptations. 
  
  \subsection{Justification de l'introduction de la quantité $R^\epsilon$}\label{SecReps}
  Le point de départ de cette étude est un calcul de \cite[p.\ 26-27]{JKO} (voir aussi \cite{BLLMilan}), qui indique que pour $u^\epsilon$, $u^*$,  $R^\epsilon$ respectivement définis par \eqref{DefUepsilon}, \eqref{DefUetoile},  et \eqref{DefResteEps}, on a
  \begin{align}
    \label{eq:Alg0}
    -\operatorname{div}\left( a\left( x/\epsilon\right)\nabla R^\epsilon(x)\right)
    &=\epsilon\, \operatorname{div} \left(a\left( \frac{x}{\epsilon} \right)\sum_{k=1}^d w_k\left( \frac{x}{\epsilon}\right) \nabla \partial_k u^*(x) \right)
    - \sum_{i=1}^d\sum_{k=1}^d M^i_k \left(\frac{x}{\epsilon}\right)  \partial_{ik} u^*(x),
  \end{align}
  où
  \begin{equation}
    \label{AssocAPot}
    M_k^i(x):=A^*_{ik} - a\left(x\right) \left( \delta_{ik}+ \partial_i w_k\left(x\right) \right).
  \end{equation}
  Pour tout $k \in [\![1,d]\!]$, comme $\operatorname{div}(M_k)=0$, il existe un potentiel $B_k =
  \left[B_k^{ij}\right]_{1\leq i,j\leq d}$ (voir \cite[p.\
  26-27]{JKO}) associé à $M_k$, c'est-à-dire une fonction antisymétrique par rapport aux indices $i$ et $j$ qui
  satisfait 
  \begin{equation}
    \label{eq:potentiel}
    \operatorname{div} \left( B_k\right)=M_k.
  \end{equation}
Grâce à \eqref{eq:Alg0}, on déduit que
  \begin{equation}
      \label{eq:Alg}
      -\operatorname{div}\left( a\left( x/\epsilon\right)\nabla R^\epsilon(x)\right)= \operatorname{div}\left(H^\epsilon(x) \right) \quad\text{dans}\quad \Omega,
    \end{equation}
    où
    \begin{align}\label{Def2Hepsilon}
      H^\epsilon_i(x)=\epsilon \sum_{k=1}^d a\left(\frac{x}{\epsilon}\right) w_k\left( \frac{x}{\epsilon}\right) \partial_{ik}u^*(x)
    -\epsilon \sum_{j,k=1}^d  B^{ij}_k\left(\frac{x}{\epsilon}\right)\partial_{jk} u^*(x).
    \end{align}
    Comme $a_{\rm per}$ est périodique et $\widetilde{a} \in {L}^r\left(\mathbb{R}^d\right)$, le potentiel $B_k$ se
    construit, pour chaque $k \in [\![1,d]\!]$, en séparant sa partie périodique $B_{k,\rm
      per}$ et sa partie due au défaut $\widetilde{B}_k$, pour laquelle on démontre que $\nabla \widetilde{B}_k \in
    {L}^r\left(\mathbb{R}^d\right)^{d\times d}$ (par la théorie de Calder\'on-Zygmund, voir \cite[p.\
    233]{Meyer2}). Par le Théorème de Morrey \cite[Th.~9.12]{brezis}, cela implique alors que, pour
    tout $i,j,k\in [\![1,d]\!],$
    \begin{align}
        \label{PotSlin}
        \left|B_k^{ij}(x)-B_k^{ij}(y)\right| \leq C |x-y|^{1-{\nu_r}}, && \forall x, y \in \mathbb{R}^d,
    \end{align}
pour $\nu_r$ défini par \eqref{DefNur}.
    Grâce à \eqref{HypoBLL}, \eqref{CorrSlin} et à \eqref{PotSlin}, on déduit de \eqref{Def2Hepsilon} que, pour
    tout $p \in [1,+\infty]$, si $f\in L^p(\Omega)$, 
    \begin{align}
      \label{HLp}
      \left\|H^\epsilon\right\|_{{L}^p(\Omega)^d} \leq C \epsilon^{{\nu_r}} \left\| \nabla^2
      u^*\right\|_{{L}^p(\Omega)^{d\times d}} \leq C \epsilon^{{\nu_r}} \left\| f \right\|_{{L}^p(\Omega)}\, .
    \end{align}
Comme on sait par ailleurs par régularité elliptique que $\nabla w_j \in {L}^\infty\left(\mathbb{R}^d\right)$, si $f\in 
C^{0,\beta}$ avec $\beta\leq \alpha$, on peut démontrer, toujours à partir de \eqref{Def2Hepsilon},
\eqref{PotSlin}, \eqref{HypoBLL}, \eqref{CorrSlin}, que
    \begin{align}
      \label{MajorH}
      \left\|H^\epsilon\right\|_{{C}^{0,\beta}(\Omega)^d} 
      \leq C\epsilon^{{\nu_r}-\beta} \left\|f\right\|_{{C}^{0,\beta}(\Omega)} 
      \quad\text{et}\quad  \left\|H^\epsilon\right\|_{{L}^\infty(\Omega)^d} \leq C\epsilon^{\nu_r} \left\|f\right\|_{{C}^{0,\beta}(\Omega)}.
    \end{align}
    
    L'enjeu de la suite de la démonstration consiste à tirer parti de \eqref{eq:Alg} et du contrôle sur $H^\epsilon$, à savoir \eqref{HLp} et \eqref{MajorH}, pour borner $R^\epsilon$.

    \subsection{Convergence dans ${H}^1(\Omega_1)$}\label{sec:conv-dans-rmh1}
    
Nous esquissons dans cette section la preuve de \eqref{PremRes} et \eqref{DeuxRes}. \red{Dans cette section, nous
utilisons le principe du maximum et certaines de ses conséquences (inégalité de Harnack, estimation de De
Giorgi-Nash-Moser). C'est pourquoi elle n'est valable, telle quelle, que pour des équations scalaires, et nécessite une adaptation pour traiter le cas de systèmes.}
      Supposons momentanément que $f \in {L}^{p}\left(\mathbb{R}^d\right)$ avec $p>d$. La fonction $R^\epsilon$
      satisfait \eqref{eq:Alg} dans $\Omega$. Néanmoins, elle n'est pas nulle au bord, ce qui empêche de faire un
      simple raisonnement variationnel. On la scinde donc en deux parties $R^\epsilon=R^\epsilon_1+R^\epsilon_2$ telles que
      \begin{align}
	&-\operatorname{div}\left( a\left( x/\epsilon\right)\nabla R^\epsilon_1(x)\right)= 0 \quad\text{dans}\quad \Omega && \quad\text{et}\quad R^\epsilon_1 =-\epsilon \sum_{j=1}^d w_j\left(\cdot/\epsilon\right)\partial_j u^* \quad\text{sur}\quad \partial \Omega,
	\label{DefR1}
	\\
	&-\operatorname{div}\left( a\left( x/\epsilon\right)\nabla R^\epsilon_2(x)\right)= \operatorname{div}\left(H^\epsilon(x) \right) \quad\text{dans}\quad \Omega && \quad\text{et}\quad R^\epsilon_2 =0 \quad\text{sur}\quad \partial \Omega.
	\label{DefR2}
      \end{align}
      Grâce au principe du maximum et à \eqref{CorrSlin}, on peut estimer $R^\epsilon_1$
      \begin{align}\label{eq:1}
        \left\|R^\epsilon_1\right\|_{{L}^\infty(\Omega)} \leq
        \left\|\epsilon \sum_{j=1}^d w_j\left(\frac{\cdot}{\epsilon}\right)\partial_j u^*\right\|_{C^0(\overline\Omega)}
        \leq C \epsilon^{{\nu_r}}\left\|u^*\right\|_{{C}^{1}(\overline\Omega)},
      \end{align}
      puis, grâce à une injection de Sobolev de ${W}^{2,p}(\Omega)$ dans $C^{1,\gamma}(\Omega)$ (pour un certain $\gamma>0$) et à
      l'estimation de régularité elliptique classique \cite[Lem.\ 9.17 p.\ 242]{GT}, on obtient
      \begin{align}
	\label{MajorR1}
	  \left\|R^\epsilon_1\right\|_{{L}^\infty(\Omega)} 
	  \leq C \epsilon^{{\nu_r}}\left\|u^*\right\|_{{W}^{2,p}(\Omega)}
	  \leq C \epsilon^{{\nu_r}} \left\|f\right\|_{{L}^p(\Omega)}.
      \end{align}
On étudie maintenant $R_2^\epsilon$.
      Rappelons que, grâce à \cite[Th.\ 1.1]{GruterWidman}, si $G^\epsilon$ est la fonction de Green relative à
      l'Equation \eqref{DefUepsilon} (donc avec conditions de Dirichlet homogènes au bord), alors la fonction
      $\nabla_y G^\epsilon(x,\cdot)$ est bornée dans l'espace de Marcinkiewicz
      $\left({L}^{\frac{d}{d-1},\infty}(\Omega)\right)^d$, uniformément en $x\in\Omega$ et en $\epsilon>0$. Ainsi, en réécrivant
      \begin{align*}
        R^\epsilon_2(x)=-\int_{\Omega} \nabla_y G^\epsilon(x,y)H^\epsilon(y) {\rm{d}} y,
      \end{align*}
      et en utilisant \eqref{HLp} (rappelons que $p>d$), on obtient
      \begin{align}
	\label{MajorR2}
        \left\|R^\epsilon_2\right\|_{{L}^\infty\left(\Omega\right)} \leq C \left\|H^\epsilon\right\|_{{L}^p(\Omega)} \leq C \epsilon^{{\nu_r}} \left\|f\right\|_{{L}^p(\Omega)},
      \end{align}
      Par conséquent, \eqref{MajorR1}, \eqref{MajorR2} et la deuxième inégalité de \eqref{eq:1} impliquent
      \begin{align}
	\label{ResultatLinfty}
	\left\|u^\epsilon-u^*\right\|_{{L}^\infty\left(\Omega \right)}
	\leq&  C\epsilon^{{\nu_r}} \left\|f\right\|_{{L}^p(\Omega)}.
      \end{align}
Ceci étant vrai pour tout $f\in L^p(\Omega)$, un argument de dualité (voir \cite[Th.\ 3.3]{KLSGreenNeumann} dans le
cas $\nu_r=1$) permet d'estimer $G^\epsilon - G^*$, où $G^*$ est la fonction de Green de l'Equation
\eqref{DefUetoile}, sur le domaine $\widetilde{\Omega}(x,y):=\Omega \cap B(y,|x-y|/16)$. On obtient ainsi
      \begin{align*}
        \left\| G^\epsilon(x,\cdot)-G^*(x,\cdot) \right\|_{{L}^{p'}\left(\widetilde{\Omega}(x,y)\right)} \leq C \epsilon^{{\nu_r}} |x-y|^{2-\frac{d}{p}-{\nu_r}} && \forall x, y \in \Omega,
      \end{align*}
      On applique alors la preuve de \cite[Lem.\ 3.2]{KLSGreenNeumann}, qui démontre une version au bord du résultat \eqref{ResultatLinfty}, que l'on applique à la fonction $G^\epsilon$
      \begin{align*}
	  \left|G^\epsilon(x,y)-G^*(x,y)\right| 
	  \leq &C |x-y|^{-d/p'} \left\|G^\epsilon(x,.)-G^*(x,.)\right\|_{{L}^{p'}\left(\widetilde{\Omega}(x,y)\right)} 
	  \\
	  &+ C \epsilon^{\nu_r} |x-y|^{1-{\nu_r}} \left\|\nabla_y G^*(x,.) \right\|_{{L}^\infty\left(\widetilde{\Omega}(x,y)\right)}
	  \\
	  &+ C \epsilon^{\nu_r} |x-y|^{2-\frac{d}{p}-{\nu_r}} \left\|\left(\nabla_y\right)^2 G^*(x,.) \right\|_{{L}^p\left(\widetilde{\Omega}(x,y)\right)}.
      \end{align*}
      Le résultat de \cite[Th.\ 1]{Dolzmann} permet de borner $\nabla_y G^*$ et $(\nabla_y)^2 G_*$, point par point. 
      Ainsi, grâce à l'inégalité de Hölder, on obtient
      \begin{align}\label{DefGreen1}
	\left|G^\epsilon(x,y)-G^*(x,y)\right| 
	\leq C  \epsilon^{{\nu_r}} |x-y|^{2-d-{\nu_r}} && \forall x, y \in \Omega.
      \end{align}
      Si on suppose maintenant, conformément à l'hypothèse du Théorème~\ref{ThDefautLp}, que $f\in {L}^2(\Omega)$, l'inégalité \eqref{DefGreen1} implique, via
      l'inégalité de Young et le fait que $|x|^{2-d-{\nu_r}}\in L^1_{\rm loc}(\mathbb R^d)$,
      \begin{equation}
        \left\|u^\epsilon-u^*\right\|_{{L}^2(\Omega)} \leq C \epsilon^{\nu_r} \left\|f\right\|_{{L}^2(\Omega)},
      \end{equation} 
      d'où \eqref{PremRes}, grâce à \eqref{CorrSlin} et au fait que $u^*\in H^1(\Omega)$.
      Cette estimation \eqref{PremRes} de la norme $L^2$ de $R_\epsilon$ se transmet en l'estimée \eqref{DeuxRes}
      de son gradient en utilisant l'inégalité de Caccioppoli et des arguments similaires à ceux ci-dessus.

    \subsection{Estimation ${L}^\infty$ sur le gradient: le cas homogène}
      
      Nous adaptons la preuve du résultat \cite[Lem.\ 16]{AvellanedaLin} qui implique que, si $-\operatorname{div}\left(a(\cdot/\epsilon)\nabla u^\epsilon\right)=0$ dans $B(0,2)$, alors 
      \begin{align}
	\label{ThAL}
	\left\|\nabla u^\epsilon\right\|_{{L}^\infty(B(0,1))} \leq C\left\|u^\epsilon \right\|_{{L}^\infty(B(0,2))}.
      \end{align}
Dans \cite{AvellanedaLin}, l'ingrédient essentiel de la preuve est le caractère borné du correcteur, propriété
impliquée par la périodicité. Cependant, l'uniforme H-convergence et la sous-linéarité du correcteur sont en fait
suffisants pour appliquer leur preuve. La démonstration de \eqref{ThAL} se fait en trois étapes: 
      \begin{enumerate}
        \item{Initialisation (voir \cite[Lem.\ 14]{AvellanedaLin}, avec un second membre nul): on obtient l'estimation 
        \begin{align}
	  &\sup_{|x| \leq \theta} \left| u^\epsilon(x) - u^\epsilon(0)- \sum_{i=1}^d 
	  \left\{x_i + \epsilon w_i\left( \frac{x}{\epsilon} \right) \right\} \fint_{B(0,\theta)} \partial_i u^\epsilon \right|
	  \leq \theta^{1+\gamma} \left( \fint_{B(0,2)} \left|u^\epsilon\right|^2 \right)^{1/2},
	  \label{EqLemIter1prime}
	\end{align}
	uniforme en $\epsilon$ suffisamment petit, à une échelle $\theta \in ]0,1[$ fixée.
	Cette étape repose sur l'existence des correcteurs et sur une propriété d'uniforme H-convergence de
    $a\left(\cdot/\epsilon\right)$ vers $A^*$, c'est-à-dire  \eqref{eq:unifHcv}.}
        \item{Itération (voir \cite[Lem.\ 15]{AvellanedaLin}): 
        on répète l'étape précédente sur les boules $B(0,\theta^2)$, $B(0,\theta^3)$, etc., jusqu'à l'échelle $\theta^k$ d'ordre $\epsilon$. Cette étape se déroule dans notre cas comme dans le cas périodique.}
        \item{\textit{Blow-up} (voir \cite[Lem.\ 16]{AvellanedaLin}): au cours de cette étape, 
            on utilise la théorie classique de Schauder pour estimer $\nabla u^\epsilon(0)$. (Ici, le point $0$ ne
            joue pas de rôle particulier.) Il faut pour cela assurer un contrôle en norme ${L}^\infty$ sur les
            correcteurs rescalés $\epsilon w_j\left(\cdot/\epsilon\right)$. De nouveau, ce contrôle est immédiat dans le cas périodique ou si $r<d$, parce que les correcteurs sont alors bornés; dans le cas où $r>d$, il se fait via la propriété de sous-linéarité sur les correcteurs $w_i$ présente dans \eqref{CorrSlin}.}
      \end{enumerate}

      \begin{remarque}\label{rk1}
      La preuve ci-dessus ne dépend pas du fait qu'on traite un système ou une équation. Ce n'est pas le
      cas de la section~\ref{sec:conv-dans-rmh1}, où nous avons utilisé le principe du maximum et des estimations
      de de Giorgi-Nash. En conséquence, pour la preuve dans le cas d'un système, il serait nécessaire d'inverser
      l'ordre des arguments : d'abord démontrer \eqref{ThAL}, puis démontrer les résultats de la
      section~\ref{sec:conv-dans-rmh1}, en utilisant, en lieu et place des estimations de de Giorgi-Nash et du
      principe du maximum, l'estimation \eqref{ThAL}.         
      \end{remarque}

    \subsection{Estimation sur les fonctions de Green}
      L'estimation suivante est démontrée dans \cite[Th.\ 1.1]{GruterWidman}, sous des hypothèses beaucoup plus
      générales que celles supposées ici :
      \begin{align}
	\label{Green1Dir}
	&\left|G^\epsilon(x,y)\right| \leq C |x-y|^{2-d} && \forall x \neq y \in \Omega.
      \end{align}
D'autre part, \eqref{ThAL}, après changement d'échelle, implique que $\|\nabla
G^\epsilon(\cdot,y)\|_{L^\infty(B_{R})}\leq C R^{-1}\|G^\epsilon(\cdot,y)\|_{L^\infty(B_{2R})},$ pour toutes boules
concentriques $B_R$ et $B_{2R}$ telles que $y\notin B_{2R}$. Ceci permet de démontrer, en utilisant $R=|x-y|/4$,
      \begin{align}
        \label{Green2Dir}
	  & \left| \nabla_x G^\epsilon(x,y) \right| \leq C |x-y|^{1-d} \quad\text{et}\quad  \left| \nabla_y G^\epsilon(x,y) \right| \leq C |x-y|^{1-d}
	  && \forall x\neq y \in \Omega_1,
      \end{align}
      et, en remarquant que, en dérivant par rapport à $y$ l'équation $-\operatorname{div}_x(a(x)\nabla_x G(x,y) =
      \delta(x-y),$ on a $-\operatorname{div}_x\left(a(x/\epsilon) \nabla_y\nabla_x G^\epsilon(x,y_0)\right)=0$ dans $\Omega
      \backslash B(y_0,R)$, pour tous $y_0 \in \Omega$ et $R>0$. On obtient donc, en utilisant \eqref{ThAL} et \eqref{Green2Dir}
      \begin{align}
	\label{Green3Dir}
	& \left| \nabla_x \nabla_y G^\epsilon(x,y) \right| \leq C |x-y|^{-d} && \forall x \neq y \in \Omega_1.
      \end{align}
      
      \subsection{Estimation ${L}^\infty$ sur le gradient: le cas non-homogène}

	On suppose désormais que $f\in {C}^{0,\beta}(\overline\Omega)$, avec $\beta\leq \alpha.$ 
	Nous appliquons alors la preuve du résultat \cite[Lem.\ 3.5]{KLSGreenNeumann}, dont nous esquissons la démonstration.
	Si $G^\epsilon$ satisfaisait \eqref{Green3Dir} sur tout le domaine, et si $R^\epsilon$ était nul sur $\partial \Omega$ (ce qui n'est pas le cas en général), alors on déduirait de \eqref{eq:Alg} que
	\begin{align}
	  \nabla R^\epsilon(x)= -\int_{\Omega} \nabla_x \nabla_y G^\epsilon(x,y)  H^\epsilon(y) {\rm{d}} y,
	\end{align}
	d'où, en isolant la singularité de $\nabla_x \nabla_y G^\epsilon$ et en utilisant \eqref{HLp} (pour $p=+\infty$) et \eqref{MajorH},
	\begin{align}
	  \label{ResLem35}
          \left\|\nabla R^\epsilon\right\|_{{L}^\infty(\Omega)} \leq&
          C \Big[ \epsilon^2 \left\|H^\epsilon\right\|_{{C}^{0,\beta}(\Omega)} + \ln \left(2+ \epsilon^{-1}\right) \left\|H^\epsilon\right\|_{{L}^\infty(\Omega)} \Big]
          \leq C \epsilon^{{\nu_r}}\ln\left(2+\epsilon^{-1}\right) \left\|f\right\|_{{C}^{0,\beta}(\Omega)}.
	\end{align}
	Pour prendre en compte que \eqref{Green3Dir} est une estimée intérieure et que $R^\epsilon$ est non nul au bord de $\partial \Omega$, on procède à une localisation en multipliant $R^\epsilon$ par une fonction de cut-off, et on démontre ainsi \eqref{TroisRes}.

\section*{Remerciements}

Le travail du second auteur a été partiellement financé par l'ONR, Grant N00014-15-1-2777, et par l'EOARD, Grant FA-9550-17-1-0294.

  \bibliographystyle{plain}
  \bibliography{Bib}

\end{document}